\documentclass{EMSpublic}

\makeatletter
\def\ps@titlepage{\let\@mkboth\@gobbletwo
 \def\@oddfoot{\hfill {\small\rm\thepage}}
 \def\@oddhead{}
 \def\@evenhead{}\let\@evenfoot\@oddfoot
 \def\sectionmark##1{}\def\subsectionmark##1{}}
\def\@maketitle{%
  \newpage
  \null
  \vskip 1pc%
  \begin{center}%
  \let \footnote \thanks
    {\large\@title \par}%
    \vskip 2pc plus 2pt minus 1pt%
    \normalsize
    \@author \par
    \vskip 4pt plus 1pt%
    {\itshape \@affiliation \par}%
    \vskip 6pt plus 1pt%
    \hphantom{(Received \@received)}\par%
    \vskip 2pc plus 2pt minus 1pt%
    \end{center}%
  \par
}
\makeatother

\usepackage{bbm}

\newcommand\C{\mathrm{C}}
\newcommand\CI{\C_1}
\newcommand\ccint[1]{\mathopen[ #1 \mathclose]}
\newcommand\ocint[1]{\mathopen] #1 \mathclose]}
\newcommand\ooint[1]{\mathopen] #1 \mathclose[}
\newcommand\cst{c}
\newcommand\D{D}
\newcommand\dd[1]{\,{\rd}#1}
\newcommand\eps\varepsilon
\newcommand\essinf{\operatorname*{ess\,inf}}
\newcommand\esssup{\operatorname*{ess\,sup}}
\newcommand\HT[1]{\mathcal{H}_{#1}}
\newcommand\Id{\mathbbm{1}}
\newcommand\ind[1]{1^{}_{#1}}
\newcommand\LL[1]{\mathrm{L}^{#1}}
\newcommand\maxp{\operatornamewithlimits{max\vphantom{p}}}
\newcommand\NN{\mathbb{N}}
\newcommand\RR{\mathbb{R}}
\newcommand\RRnn{\RR_{\ge0}}
\newcommand\tocc{\xrightarrow{\textup{cc}}}
\newcommand\umax{u_{\mathrm{max}}}
\newcommand\wmin{w_{\mathrm{min}}}
\newcommand\wo\backslash

\newcommand\es\,  % space before punctuation marks in numbered equations
\newcommand\mt[1]{\boldsymbol{#1}}  % matrix
\def\nrm#1{{\hbox{$\left#1\vbox to7pt{}\right.\nulldelimiterspace0pt$}}}
\newcommand\scal[1]{\langle #1 \rangle}
\newcommand\st[1]{\mathcal{#1}}     % set
\newlength\thnormheight
\newcommand\thnorm[1]{%
  \mathopen{\mkern-\thickmuskip}\mathopen\nrm\bracevert\mkern-\thickmuskip #1
  \settoheight\thnormheight{A}%
  \rule{0ex}{\thnormheight}%
    \mkern-\thickmuskip \mathclose\nrm\bracevert\mathclose{\mkern-\thickmuskip}}
\newcommand\tnorm[1]{|\!|\!|#1|\!|\!|}
\newcommand\trp{^\mathrm{T}}
\newcommand\vc[1]{\boldsymbol{#1}}  % vector

\newcommand\pr[1]{\textup{(}#1\textup{)}}
\newcommand\textdef[1]{\emph{#1}}

\newcommand\dref[1]{Definition~\textup{\ref{def:#1}}}

\newcommand\eref[1]{\textup{(\ref{eq:#1})}}

\newcommand\lref[1]{Lemma~\textup{\ref{lem:#1}}}

\newcommand\pref[1]{Proposition~\textup{\ref{prop:#1}}}

\newcommand\sref[1]{Section~\textup{\ref{sec:#1}}}

\newcommand\thref[1]{Theorem~\textup{\ref{thm:#1}}}
\newcommand\Thref[1]{Theorem~\textup{\ref{thm:#1}}}

\makeatletter
\newcommand\cf{cf.\@ifnextchar,{}{\ }}
\newcommand\eg{e.g.\@ifnextchar,{}{, }}
\newcommand\ie{i.e.\@ifnextchar,{}{, }}
\makeatother

\newcommand\chap{Ch.~}
\newcommand\etal{et al.}
\newcommand\exm{Exm.~}
\newcommand\lem{Lemma }
\newcommand\lems{Lemmas }
\newcommand\lemmas{Lemmas }
\newcommand\prop{Prop.~}
\newcommand\props{Props.~}
\newcommand\propositions{Propositions }
\newcommand\satz\thm
\newcommand\saetze\thms
\newcommand\sect{Sec.~}
\newcommand\sects{Secs.~}
\newcommand\sections{Sections }
\newcommand\Tctp{This completes the proof.}
\newcommand\thm{Thm.~}
\newcommand\thms{Thms.~}
\newcommand\wrt{with respect to }

\makeatletter
\def\itemlabel#1{\def\@currentlabel{\@itemlabel}\label{#1}}
\newcommand\itlabel[2]{\textup{#1}%
  \def\@currentlabel{\textup{#1}}%
  \label{#2}}
\makeatother

\makeatletter
\renewcommand\thesubsection{\thesection.\@arabic\c@subsection}
\renewcommand\subsection{\renewcommand{\@seccntformat}[1]{%
   \csname thesubsection\endcsname.\hspace{0.5em}}%
   \@startsection{subsection}{2}{\z@}{3.25ex\@plus 1ex \@minus .2ex}%
                                     {1.4ex \@plus .2ex}%
                                     {\normalfont\bfseries%
                                      \mathversion{bold}}}
\def\@citex[#1]#2{%
  \let\@citea\@empty
  \@cite{\@for\@citeb:=#2\do
    {\@citea\def\@citea{,\penalty\@m\ }%
     \edef\@citeb{\expandafter\@firstofone\@citeb\@empty}%
     \if@filesw\immediate\write\@auxout{\string\citation{\@citeb}}\fi
     \@ifundefined{b@\@citeb}{\mbox{\reset@font\bfseries ?}%
       \G@refundefinedtrue
       \@latex@warning
         {Citation `\@citeb' on page \thepage \space undefined}}%
       {\hbox{\bfseries\csname b@\@citeb\endcsname}}}}{#1}}
\makeatother
\newtheorem*{ack}{Acknowledgements}

\begin{document}

\title[Discrete Approximation of Non-Compact Operators]{Discrete Approximation
  of Non-Compact Operators Describing Continuum-of-Alleles Models} 
\author{Oliver Redner}
\affiliation{Institut f\"ur Mathematik und Informatik,
  Universit\"at Greifswald, Jahnstr.\ 15a, 17487~Greifswald, Germany
  \textup{(redner@uni-greifswald.de)}}
\maketitle

\begin{abstract}
  We consider the eigenvalue equation for the largest eigenvalue of certain
  kinds of non-compact linear operators given as the sum of a multiplication
  and a kernel operator.  It is shown that, under moderate conditions, such
  operators can be approximated arbitrarily well by operators of finite rank,
  which constitutes a discretization procedure.  For this purpose, two
  standard methods of approximation theory, the Nystr\"om and the Galerkin
  method, are generalized.  The operators considered describe models for
  mutation and selection of an infinitely large population of individuals that
  are labeled by real numbers, commonly called continuum-of-alleles (COA)
  models.
\end{abstract}


% % Eigenvalue problems
    {47B34}  % Kernel operators

\section{Introduction}

This article is concerned with eigenvalue equations on $\LL{1}(I)$ of the form
\begin{equation}
  \label{eq:coaequil}
  r(x) \, p(x) + \int_I \bigl[ u(x,y) \, p(y) - u(y,x) \, p(x) \bigr] \dd{y} =
  \lambda \, p(x) 
  \qquad\text{for all $x \in I$.}
\end{equation}
Here, $p$ is a probability density on the set $I$, which is either taken to be
a compact interval $\ccint{a,b}$ or the real line $\RR$, \ie $p \in \LL{1}(I)$
with $p\ge0$ and $\int_I p(x) \dd{x} = 1$.  Sufficient conditions for the
existence and uniqueness of solutions of \eref{coaequil} were given by
B\"urger, see \cite[\chap IV.3]{Bue00}, in which case $\lambda$ is the largest
eigenvalue.

If one is interested in a discrete approximation of \eref{coaequil}, one faces
the problem that the operator acting on $p$ is the sum of a multiplication
operator and a kernel operator; and the former is never compact (apart from
trivial cases).  Therefore, a direct application of most standard methods of
approximation theory fails because, for these, compactness is a prerequisite.
In this article, it will be shown that, under some moderate extra conditions,
these methods can nevertheless be applied.

One motivation to study equations of the form \eref{coaequil} is their
occurrence in population genetics, which is concerned with the (micro)evolution
of the genetic composition of populations.  For many situations, individuals
are adequately described by a continuous scalar variable, representing, for
example, a quantitative character under selection.  This leads to the
definition of so-called continuum-of-alleles (COA) models, in which
individuals are identified with this variable, referred to as their type.
Usually, selection is then modeled by type-dependent fitness values, whereas
mutation is described, for every source type, by a probability distribution
for the mutant types.  For a recent review of and relevant literature on COA
models, see \cite{Bue00}.

In population genetics, evolution may quite generally be assumed to proceed in
continuous time, with overlapping generations, or in discrete generations.
For the COA model, in both cases, equilibrium is described by an equation of
the form \eref{coaequil}, compare \cite{Kim}.  Here, $I$ is the set of
possible types.  Assuming the population to be effectively infinite, we
represent it by the probability density $p$.

The notation chosen here best fits the case of continuous time, where
$r(x)$ describes the effective reproduction rate of type $x$ (\ie the
difference of its birth and death rate), the so-called Malthusian
fitness, and $u(x,y)$ is the mutation rate $u_1(y)$ of type $y$ times
the density $m(x,y)$ of mutant types $x$, conditioned on a mutation to
occur for $y$.  With discrete generations, $r(x)$ has the
interpretation of the expected number of offspring of an individual of
type $x$, \ie its Wrightian fitness, and mutation is assumed to occur
during reproduction with some probability $\mu(y)$ for type $y$.  The
distribution of mutant types is again given by $m(x,y)$, hence $u(x,y)
= m(x,y) \, \mu(y) \, r(y)$.  In both cases, $\lambda$ equals the
equilibrium mean fitness $\int_I r(x) \, p(x) \dd{x}$.

There are several reasons why it is desirable to approximate a COA model by a
model with discrete types.  One reason is the need for numerical
investigations of COA models, since most of them are not tractable
analytically.  These inevitably require a discrete formulation of the model.
Another reason is that recently a simple characterization of the equilibrium
of discrete mutation--selection models has been found \cite{HRWB} (see also
\cite{GeBa03,GaGr,BBBK}); this takes the form of a scalar maximum principle in
a limit of infinitely many types that densely fill a compact interval.
Gaining a better understanding of the relation between models with discrete
and continuous types is therefore promising to enable a transfer of some of
these results.

This article starts with a summary of B\"urger's results on \eref{coaequil} in
\sref{coageneral}, since these form the basis for our treatment.  We will then
consider two methods to approximate compact kernel operators and extend them
to our case.  One, the Nystr\"om method, is applicable to continuous functions
$r$ and $u$ on compact intervals $I$ and involves sampling\footnote{The term
  \textdef{sampling} is used in the meaning also used in signal processing:
  Instead of a continuous function one considers its values at a (properly
  chosen) finite set of points.}  of these functions.  This is presented in
\sref{compact}.  The other one, the Galerkin method, is based on projections
to finite-dimensional subspaces and works---in principle---for a broad class
of compact operators.  In our case, however, one has to make relatively strong
assumptions, \eg that the functions $r$ and $u$ are, in some sense, uniformly
continuous.  Then, it turns out that the local averaging in the projection
process can be replaced by sampling again (if an additional condition is
satisfied).  This is discussed in \sref{unbounded}.  A comparison of both
methods in \sref{comparison} and an outlook in \sref{outlook} complete this
article.

\section{General properties}
\label{sec:coageneral}
\label{sec:operatornot}

Let us first put the equilibrium condition \eref{coaequil} in operator
notation.  Since we are interested in probability densities, we will
consider $\LL{1}(I)$, or a subspace thereof, as the underlying
function space.  We define the total mutation rate of type $x$ as
\begin{equation}
  \label{eq:u1}
  u_1(x) = \int_I u(y,x) \dd{y}
\end{equation}
and, for notational brevity,
\begin{displaymath}
  w = u_1 - r \es.
\end{displaymath}
Then, \eref{coaequil} is equivalent to the eigenvalue equation
\begin{equation}
  \label{eq:coaequilop}
  (A + \lambda) p = 0 \es,
\end{equation}
where, for elements $f$ of the function space and all $x \in I$,
\begin{align}
  \label{eq:deft}
  (Tf)(x) &= w(x) f(x) \es, \\
  \label{eq:defu}
  (Uf)(x) &= \int_I u(x,y) f(y) \dd{y} \es, \\
  \label{eq:defa}
  A &= T-U \es.
\end{align}
As mentioned above, being a (non-zero) multiplication operator, $T$
cannot be compact (compare \cite[\thm 2.1]{Mor}).  Strong results like
analogs to the Perron--Frobenius theorem, however, are only available
for compact, or at least power compact\footnote{An operator is said to
  be \textdef{power compact} if one of its powers is compact.},
operators, see Schaefer \cite[\chap V]{Schae}.  Therefore one
considers the following family of kernel operators:
\begin{displaymath}
  (K_\alpha f)(x) = \int_I k_\alpha(x,y) f(y) \dd{y} \es,
\end{displaymath}
where
\begin{displaymath}
  k_\alpha(x,y) = \frac{u(x,y)}{w(y)+\alpha} \es.
\end{displaymath}
These are, under conditions that will be given shortly, power compact
or even compact.  Their connection to the operator $A$ from
\eref{defa} is stated in the following
\begin{lemma}\textup{\cite[\prop 2.1(i)]{Bue88}}
  \label{lem:connectka}
  Let\/ $T$, $U$ be operators in a Banach space\/ $X$, with\/ $U$
  being bounded, $T$ densely defined, \ie $\overline{\D(T)} = X$,
  and\/ $T+\alpha$ invertible.  Then\/ $f$ is an eigenvector of\/ $A =
  T-U$ with eigenvalue\/ $-\alpha$, \ie $0 \neq f \in \D(A) = \D(T)$
  and
  \begin{displaymath}
    (A+\alpha)f=0 \es,
  \end{displaymath}
  if and only if\/ $g = (T+\alpha) f$ is an eigenvector of\/ $K_\alpha
  = U(T+\alpha)^{-1}$ with eigenvalue\/ $1$,
  \begin{displaymath}
    (K_\alpha-1)g=0 \es.
  \end{displaymath}
\end{lemma}
So, explicitly in our case, the eigenvalue equation \eref{coaequilop}
is equivalent to
\begin{equation}
  \label{eq:coaequilopk}
  (K_\lambda - 1) q = 0
\end{equation}
with $q = (T+\lambda) p$.  This equation can now be used to find
sufficient conditions for the existence and uniqueness of a solution
of \eref{coaequilop}.

An important class of bounded kernel operators from $\LL{q}(I)$ into
$\LL{p}(I)$ ($1 \le p$, $q \le \infty$) are the Hille--Tamarkin
operators, see \cite[\sect 11.3]{Joe}.  Their kernels need to satisfy
\begin{equation}
  \label{eq:hiltampq}
  \thnorm{K}^{}_{pq} := \|k_1\|^{}_p < \infty
  \qquad\text{with}\qquad
  k_1(x) = \|k(x,.)\|^{}_{q'} \es,
\end{equation}
where $(Kf)(x) = \int_I k(x,y) \, f(y) \dd{y}$, $k(x,.)$ denotes the
function $y \mapsto k(x,y)$, and $q'$ is the conjugate exponent to $q$
satisfying $\frac1q + \frac1{q'} = 1$, $1 \le q' \le \infty$.  The
Hille--Tamarkin norm $\thnorm{.}^{}_{pq}$ turns the set $\HT{pq}(I)$
of all Hille--Tamarkin operators into a Banach space \cite[\satz
11.5]{Joe}.  Here, we are interested in $p=q=1$, in which case
\eref{hiltampq} yields
\begin{displaymath}
  \thnorm{K}^{}_{11} = 
  \int_I \esssup_{y \in I} |k(x,y)| \dd{x} < \infty
\end{displaymath}
and $K^2$ is compact for every $K \in \HT{11}(I)$ \cite[\satz
11.9]{Joe}.

Let us now turn to kernel operators that are power compact, positive,
and irreducible.  An operator is called \textdef{positive} if it maps
the set of non-negative functions into itself, for which, in the case
of kernel operators, non-negativity of the kernel is necessary and
sufficient \cite[p.\ 122]{Joe}.  A kernel operator is
\textdef{irreducible} if its kernel satisfies \cite[\exm 4 in \sect
V.6]{Schae}
\begin{displaymath}
  \int_{I \wo J} \int_J k(x,y) \dd{x} \dd{y} > 0
  \qquad\text{for all measurable $J \subset I$ with $|J|$, $|I \wo J| > 0$.}
\end{displaymath}
Here, $|J|$ denotes the Lebesgue measure of a measurable set $J$.
Then, the theorem of Jentzsch \cite[\thm V.6.6]{Schae}, which parallels
the Perron--Frobenius theorem for matrices, states that the spectral
radius is an algebraically simple eigenvalue with an (up to
normalization) unique positive eigenfunction (\ie strictly positive
a.e.\footnote{The abbreviation `a.e.'\ stands for `almost every' or
  `almost everywhere' and means that the set at which the condition it
  refers to is not fulfilled has zero (Lebesgue) measure.})\ and the
only eigenvalue with a positive eigenfunction.

In our case, the following requirements are sufficient for the
$K_\alpha$ to be Hille--Tamarkin operators \cite[\sect 3]{Bue88}.
\begin{itemize}
\item[\itlabel{(U1)}{it:u1}] $u$ is non-negative and measurable.
\item[\itlabel{(U2)}{it:u2}] $u_1(x)$ from \eref{u1} exists for a.e.\ $x \in
  \RR$ and $u_1 \in \LL{\infty}(I)$, \ie $u_1$ is essentially bounded.  (By
  H\"older's inequality, this implies that $U$ is bounded, \cf \cite[\prop
  3.1(ii)]{Bue88}.)
\item[\itlabel{(T1)}{it:t1}] $w = u_1-r$ is measurable and satisfies
  $\essinf_{x \in I} w(x) = 0$.  (The latter can be achieved, without
  loss of generality, by adding a suitable constant to $r$.)
\item[\itlabel{(T2)}{it:t2}] $(w+1)^{-1} \in \LL{\infty}(I)$ is then
  already a consequence of \ref{it:t1}.
\item[\itlabel{(U4)}{it:u4}] $\int_I \esssup_{y \in I}
  u(x,y)/(w(y)+\alpha) \dd{x} < \infty$ for one (and then for all)
  $\alpha>0$.
\end{itemize}
For $\alpha>0$, $K_\alpha$ is irreducible if $U$ is \cite[proof of \thm
2.2(c)]{Bue88}, \ie
\begin{equation}
  \label{eq:irreduc}
  \int_{I \wo J} \int_J u(x,y) \dd{x} \dd{y} > 0
  \qquad\text{for all measurable $J \subset I$ with $|J|$, $|I \wo J| > 0$.}
\end{equation}
To keep the equilibrium distribution from having atoms, we assume that
there is a set $J \subset I$ with positive measure for which
$\essinf_{x \in J} w(x) = 0$ such that
\begin{equation}
  \label{eq:cusp}
  \essinf_{x,y \in J} u(x,y) \int_J (w(x))^{-1} \dd{x} > 1
\end{equation}
or the integral diverges \cite[cond.\ 3'' in \sect IV.3]{Bue00}.

Putting everything together, we have the following
\begin{theorem}[B\"urger]
  \label{thm:existunique}
  Under the above conditions, \eref{coaequil} has a unique positive
  solution\/ $p \in \LL{1}(I)$ with\/ $\|p\|^{}_1 = 1$, for which\/
  $\lambda>0$ is the largest spectral value of\/ $-A$ from
  \eref{defa}.
\end{theorem}
\begin{proof}
  See the above, \cite[\thm 3.5]{Bue88}, and \cite[\sect IV.3]{Bue00}.
\end{proof}
Note that, due to \ref{it:t1}, $p$ is positive if and only if $q = (w+\alpha)
p$ is, for $\alpha > 0$.

Another result that will be needed in the sequel is
\begin{lemma}\textup{\cite[\lems 1--3 and \thm 2.2(ii)]{Bue88}}
  \label{lem:alphalambda}
  Under the above conditions, the spectral radius\/ $\rho(K_\alpha)$
  is, as a function of\/ $\alpha$, strictly decreasing and satisfies\/
  $\rho(K_\lambda) = 1$ as well as\/ $\lim_{\alpha\to\infty}
  \rho(K_\alpha) = 0$.  Thus, $\rho(K_\alpha) < 1$ implies\/ $\alpha >
  \lambda$ and\/ $\rho(K_\alpha) > 1$ implies\/ $\alpha < \lambda$.
\end{lemma}

Throughout the rest of this article, all the above criteria are
assumed to be satisfied, namely \ref{it:u1}, \ref{it:u2}, \ref{it:u4},
\ref{it:t1}, \ref{it:t2}, \eref{irreduc}, and \eref{cusp}.

\section{Discretization---compact interval}
\label{sec:compact}

Let the interval $I$ be compact and $\C(I)$ denote the Banach space of
bounded, continuous functions equipped with the supremum norm
$\|f\|^{}_\infty = \sup_{x \in I} |f(x)|$.  We consider operators $K$
of the form
\begin{equation}
  \label{eq:kernelopi}
  (Kf)(x) = \int_I k(x,y) f(y) \dd{y}
  \qquad\text{for all $x \in I$}
\end{equation}
with a continuous kernel $k\colon I \times I \to \RR$.  First
note these two basic results:
\begin{proposition}
  \label{prop:klinc}
  Any\/ $K$ of the form \eref{kernelopi} maps\/ $\LL{1}(I)$ into\/
  $\C(I) \subset \LL{1}(I)$.
\end{proposition}
\begin{proof}
  We follow the proof of \cite[\satz 2.1]{Eng}, where this is shown
  for $\LL{2}(I)$, which, since $I$ is compact, is a subspace of
  $\LL{1}(I)$.  Let $f \in \LL{1}(I)$ and $x$, $\xi \in I$ be given.
  Then,
  \begin{displaymath}
    |(Kf)(x)-(Kf)(\xi)| \le
    \int_I |k(x,y)-k(\xi,y)| \, |f(y)| \dd{y} \le
    \sup_{y \in I} |k(x,y)-k(\xi,y)| \, \|f\|^{}_1 \es.
  \end{displaymath}
  Due to the uniform continuity of $k$ in $I \times I$, we have
  \begin{displaymath}
    \operatornamewithlimits{lim\vphantom{p}}_{\xi \to x} 
    \sup_{y \in I} |k(x,y)-k(\xi,y)| = 0 \es,
  \end{displaymath}
  from which the continuity of $Kf$ follows.
\end{proof}
\begin{proposition}
  \label{prop:kcompact}
  An operator\/ $K$ of the form \eref{kernelopi} is compact from\/
  $\C(I)$ or\/ $\LL{1}(I)$ to either of the two spaces.
\end{proposition}
\begin{proof}
  Follow the proof of \cite[\satz 2.10]{Eng} (or \cite[XVII.4]{Lan}),
  where this is shown for $\LL{2}(I) \subset \LL{1}(I)$, and use
  H\"older's inequality whenever the Cauchy--Schwarz inequality is
  used.  Alternatively, see \cite[\exm 3 in \sect IV.10]{Schae}.
\end{proof}

Thus, if in our case the functions $r$ and $u$ are continuous, also
the kernel $k_\alpha$ is, for every $\alpha>0$.  It then follows from
\pref{klinc} that the equilibrium density $p$ is continuous as well.
Therefore we can restrict our attention to $\C(I)$ in our quest for a
solution of the eigenvalue equation \eref{coaequilop}.  This makes the
Nystr\"om method applicable as a discretization procedure, which will
be presented now.

\subsection{The Nystr\"om method}

The Nystr\"om method is based on quadratures, which are used for
numerical integration, \cf, \eg, Kress \cite[\chap 12]{KreLIE}.  We
will use this (slightly restricted)
\begin{definition}
  \label{def:quad}
  A \textdef{quadrature rule} $Q_n$ is a
  mapping of the form
  \begin{displaymath}
    Q_n\colon \C(I) \to \RR \es, \qquad 
    f \mapsto Q_n f = \sum_{k=1}^{N_n} \alpha_{n,k} f(t_{n,k}) \es,
  \end{displaymath}
  with $n\in\NN$, $N_n\in\NN$, \textdef{quadrature points} $t_{n,k}
  \in I$, and \textdef{quadrature weights} $\alpha_{n,k} > 0$, for $k
  \in \st{N}_n := \{1,\ldots,N_n\}$.  A sequence of quadrature rules,
  or simply a \textdef{quadrature}, $(Q_n)$ is said to be
  \textdef{convergent} if
  \begin{equation}
    \label{eq:quadprop}
    Q_n f \to Qf \qquad\text{for all $f \in \C(I)$,}
  \end{equation}
  where $Q \colon \C(I) \to \RR$ is the linear functional that
  assigns to each $f$ its integral, \ie $Qf = \int_I f(x) \dd{x}$.
\end{definition}
Another notion that is important for the Nystr\"om method is the
collectively compact convergence of operators.  The standard reference
for this matter is \cite{Ans}.
\begin{definition}
  \label{def:collcomp}
  A sequence $(K_n)$ of \pr{compact} operators in a Banach space $X$ is
  \textdef{collectively compact} if the set $\{K_n B : n\in\NN\}$ is
  relatively compact \pr{\ie its closure is compact} for every bounded set $B
  \subset X$.  If furthermore the sequence converges pointwise to an operator
  $K$ one speaks of \textdef{collectively compact convergence}, in symbols
  $K_n \tocc K$.
\end{definition}
As a direct consequence of this definition, $K$ is compact (as well as
all $K_n$).  The central result for the Nystr\"om method is
\begin{theorem}
  \label{thm:nystroem}
  Let\/ $K$ be a compact kernel operator of the form \eref{kernelopi}
  whose eigenvalue equation
  \begin{equation}
    \label{eq:ieigvaleqk}
    (K - \nu) g = 0
  \end{equation}
  is to be approximated.  To this end, let\/ $(Q_n)^{}_{n\in\NN}$ be a
  convergent quadrature with the notation as in \dref{quad}.  A
  \textdef{complete discretization} is
  given by the\/ $N_n \times N_n$ matrices\/ $\mt{K}_n$ with entries
  \begin{displaymath}
    K_{n,k\ell} = \alpha_{n,\ell} \, k(t_{n,k},t_{n,\ell}) \es,
  \end{displaymath}
  a \textdef{partial discretization} by means of the operators\/ $K_n$
  on\/ $\C(I)$ with
  \begin{equation}
    \label{eq:partdiscop}
    (K_n f)(x) = 
      \sum_{k=1}^{N_n} \alpha_{n,k} \, k(x,t_{n,k}) f(t_{n,k}) =
      Q_n(k(x,.)f) \es.
  \end{equation}
  Consider the corresponding eigenvalue equations
  \begin{equation}
    \label{eq:ieigvaleqkn}
    (\mt{K}_n - \nu_n) \vc{g}_n = 0 
    \qquad\text{and}\qquad
    (K_n - \nu_n) g_n = 0 \es,
  \end{equation}
  where\/ $\vc{g}_n$ is an\/ $N_n$-dimensional vector with
  components\/ $g_{n,k}$, and\/ $g_n \in \C(I)$.  Then, under the
  above conditions the following statements are true:
  \begin{enumerate}\renewcommand\labelenumi{\textup{(\alph{enumi})}}
  \item\itemlabel{it:ieveequiv} Both eigenvalue equations in
    \eref{ieigvaleqkn} are equivalent and connected via
    \begin{equation}
      \label{eq:partdiscvec}
      g_n(x) =
        \sum_{k=1}^{N_n} \alpha_{n,k} \, k(x,t_{n,k}) g_{n,k} \es.
    \end{equation}
  \item\itemlabel{it:ieigval} For every\/ $\nu \neq 0$ from
    \eref{ieigvaleqk} there is a sequence\/ $(\nu_n)$ of eigenvalues
    of \eref{ieigvaleqkn} such that\/ $\nu_n \to \nu$ as\/
    $n\to\infty$. Conversely, every non-zero limit point of any
    sequence\/ $(\nu_n)$ of eigenvalues of \eref{ieigvaleqkn} is an
    eigenvalue of \eref{ieigvaleqk}.
  \item\itemlabel{it:ieigvec} Every bounded sequence\/ $(g_n)$ of
    eigenfunctions of \eref{ieigvaleqkn} associated with eigenvalues\/
    $\nu_n \to \nu \neq 0$ contains a convergent subsequence; the
    limit of any convergent subsequence\/ $(g_{n_i})^{}_i$ is an
    eigenfunction of \eref{ieigvaleqk} associated with the
    eigenvalue\/ $\nu$ \pr{unless the limit is zero}.
  \end{enumerate}
\end{theorem}
\begin{proof}
  \ref{it:ieveequiv} is the statement of \cite[\thm 12.7]{KreLIE} or
  \cite[\lem 3.15]{Eng}.  \ref{it:ieigval} and \ref{it:ieigvec} rely on $K_n
  \tocc K$, which is shown, \eg, in \cite[\props 2.1, 2.2]{Ans}, \cite[\thm
  12.8]{KreLIE}, or \cite[\satz 3.22]{Eng}.  The statements then follow from
  \cite[\thms 4.11, 4.17]{Ans}.
\end{proof}

We will restrict ourselves to quadratures that allow for disjoint partitions
of $I$ with intervals $I_{n,k}$, \ie $I_{n,k} \cap I_{n,\ell} \neq \emptyset$
and $\bigcup_{k=0}^{N_n} I_{n,k} = I$, such that $t_{n,k} \in I_{n,k}$ and
$|I_{n,k}| = \alpha_{n,k}$ (with $k \in \st{N}_n$).  For such quadratures it
is easy to see that\footnote{If not noted otherwise, the following convention
  for operator norms is used.  If an operator maps a space $X$ into itself, we
  denote its norm by the same symbol as the norm of $X$, \eg $\|.\|^{}_X$, or
  $\|.\|^{}_1$ for $\LL{1}$; in all other cases the unornamented symbol
  $\|.\|$ is used.}
\begin{equation}
  \label{eq:normqn}
  \|Q_n\| = \sum_{k=1}^{N_n} \alpha_{n,k} = |I|
\end{equation}
and that the partitions are unique (up to the boundary points of the
intervals).  Furthermore we have
\begin{lemma}
  \label{lem:fine}
  Let\/ $(Q_n)$ be a convergent quadrature that allows for partitions
  of\/ $I$ as described above.  Then\/ $\lim_{n\to\infty}
  \max_{k\in\st{N}_n} |I_{n,k}| = 0$.
\end{lemma}
\begin{proof}
  Assume the contrary.  Then there are an $\eps>0$ and sequences
  $(n_i)^{}_i$ and $(k_i)^{}_i$ with $\lim_{i\to\infty} n_i = \infty$
  such that $|I_{n_i,k_i}| \ge \eps$.  Due to the compactness of $I$,
  these can be chosen in a way that $\lim_{i\to\infty} t_{n_i,k_i} =:
  t$ exists.  Now consider the continuous function $f(x) =
  \max\{1-2|x-t|/\eps,0\}$.  For this we have $Qf \le \eps/2$, but
  $\lim_{i\to\infty} Q_{n_i}f \ge \eps \lim_{i\to\infty}
  f(t_{n_i,k_i}) = \eps$, which contradicts the convergence of the
  quadrature \eref{quadprop}.
\end{proof}

\subsection{Application to the COA model}
\label{sec:iapplic}

In our case of the COA model with a compact interval $I$ and
continuous functions $r$ and $u$, the complete
discretization is given by the
following $N_n \times N_n$ matrices:
\begin{align}
  \label{eq:deftn}
  &T_{n,k\ell} = \delta_{k\ell} \, w(t_{n,k}) \ge 0 \es, \\
  \label{eq:defun}
  &U_{n,k\ell} = \alpha_{n,\ell} \, u(t_{n,k},t_{n,\ell}) \ge 0 \es, \\
  &\mt{A}_n = \mt{T}_n-\mt{U}_n, \qquad
  \mt{K}_{\alpha,n} = \mt{U}_n (\mt{T}_n+\alpha)^{-1}
  \quad\text{for $\alpha>-\min_{k\in\st{N}_n} w(t_{n,k})$.} \notag
\end{align}
The eigenvalue equations to be solved are
\begin{displaymath}
  (\mt{A}_n^{} + \lambda_n^{}) \vc{p}_n^{} = 0
  \qquad\text{with $\vc{p}_n^{} > 0$.}
\end{displaymath}
Here, $-\mt{A}_n + \cst$ is positive with a suitable constant $\cst$.  We
further have to assume that the $\mt{A}_n$ are irreducible (which
might not be the case for special choices of the $t_{n,k}$, \eg if
$u_1(t_{n,k})=0$ for some $k$).  Then, due to the Perron--Frobenius
theorem, there exist (up to normalization) unique positive
$\vc{p}_n^{}$ belonging to the eigenvalues $-\lambda_n =
-\rho(-\mt{A}_n+\cst)+\cst$, where $\rho(\mt{M})$ denotes the spectral
radius of a matrix $\mt{M}$.  With $\vc{q}_n^{} = (\mt{T}_n^{} +
\lambda_n^{}) \vc{p}_n^{}$ also the eigenvalue equations
\begin{equation}
  \label{eq:ieigvaleqqn}
  (\mt{K}_{\lambda_n,n}^{} - 1) \vc{q}_n^{} = 0
\end{equation}
are solved (and vice versa), \cf \lref{connectka}.

Both $\mt{K}_{\lambda_n,n}$ and ${\vc{q}}_n$ can be embedded into
$\C(I)$ as described by \eref{partdiscop} and \eref{partdiscvec}.
Then, with \thref{nystroem}, one might conclude the convergence
$\|q_n^{}-q\|^{}_\infty \to 0$.  In the end, however, we are
interested in the population vectors $\vc{p}_n^{}$ and their
convergence to the density $p$.  It might be easiest to interpret the
vectors $\vc{p}_n^{}$ as point measures on $I$.  But then the best one
can hope for is weak convergence since the set of point measures is
closed under the total variation norm.  It will turn out that we can
indeed achieve norm convergence if we embed the $\vc{p}_n^{}$ into
$\LL{1}(I)$ the following way.  We choose a disjoint partition of $I$
as above and let
\begin{displaymath}
  p_n^{} = \sum_{k=1}^{N_n} p_{n,k}^{} \ind{I_{n,k}} \es,
\end{displaymath}
where $\ind{J}$ denotes the characteristic function of a set $J$.
(Note that $p_{n,k}^{}$ denotes the $k$-th component of $\vc{p}_n^{}
\in \RR^{N_n}$, whereas $p_n^{}$ is an $\LL{1}$ function.)  Thus the
$\vc{p}_n^{}$ can be interpreted as probability densities on $I$, if
we normalize them such that $\|p_n^{}\|^{}_1 = 1$.  This is most
easily expressed using the induced norm $\|\vc{f}\|^{}_{(n)} :=
\sum_{k=1}^{N_n} \alpha_{n,k} |f_k|$ on $\RR^{N_n}$.  Convergence in
total variation then corresponds to $\|p_n^{}-p\|^{}_1 \to 0$
\cite[\thm 6.13]{RudRCA}.\footnote{One may also define operator analogs of the
  $\mt{A}_n$, see \cite[\sect II.2.1.2]{Diss}.}

\subsection{Convergence of eigenvalues and eigenvectors}
\label{sec:ieigvalvec}

We now come to prove the main approximation result:
\begin{theorem}
  \label{thm:ieigvalvec}
  With the notation and assumptions from \sections
  \textup{\ref{sec:operatornot}} and \textup{\ref{sec:iapplic}},
  \begin{enumerate}\renewcommand\labelenumi{\textup{(\alph{enumi})}}
  \item\itemlabel{it:thmieigval} $\lim_{n\to\infty} \lambda_n =
    \lambda > 0$ and
  \item\itemlabel{it:thmieigvec} $\lim_{n\to\infty} \|p_n^{}-p\|^{}_1
    = 0$, \ie the probability measures corresponding to these
    densities converge in total variation.
  \end{enumerate}
\end{theorem}
The idea of the proof is as follows.  In the following two lemmas, we first
determine an upper and a lower bound for the $\lambda_n$ and conclude that
there is a convergent subsequence.  Then we show that every convergent
subsequence converges to $\lambda$ and hence the sequence itself.  By
\thref{nystroem}, this implies the convergence of a subsequence of
$(q_n^{}/\|q_n^{}\|_\infty^{})$ to a (non-negative) limit function.  Since,
due to \thref{existunique}, the latter is unique, we conclude that it is
$q/\|q\|_\infty^{}$.  With this, part \ref{it:thmieigvec} can be shown.
\begin{lemma}
  \label{lem:limsup}
  There is a constant\/ $M>0$ such that\/ $|\lambda_n| \le M$ for
  all\/ $n\in\NN$.
\end{lemma}
\begin{proof}
  Using \eref{deftn} and \eref{defun}, one checks
  \begin{align*}
    |\lambda_n^{}| = 
      \frac{\|\lambda_n^{} \vc{p}_n^{}\|^{}_{(n)}}{\|\vc{p}_n^{}\|^{}_{(n)}} =
      \frac{\|\mt{A}_n^{}\vc{p}_n^{}\|^{}_{(n)}}{\|\vc{p}_n^{}\|^{}_{(n)}} &\le
      \sup_{\|\vc{f}\|^{}_{(n)}=1} \sum_{k=1}^{N_n} \alpha_{n,k}^{}
        \left| \sum_{\ell=1}^{N_n} 
          (T_{n,k\ell}^{}-U_{n,kl}^{}) f_\ell^{} \right| \\
    &\le \max_k w(t_{n,k}^{}) + \max_{k,\ell} u(t_{n,k}^{},t_{n,\ell}^{}) 
        \sum_{k=1}^{N_n} \alpha_{n,k}^{} \\
    &\le \|w\|^{}_\infty + \|u\|^{}_{\C(I \times I)} \sup_m \|Q_m^{}\| =:
      M > 0 \es.
\end{align*}
Here, $\|Q_m\| = |I|$ due to \eref{normqn}.  More generally, $\sup_m
\|Q_m\| < \infty$ holds for any convergent quadrature according to the
theorem of Banach--Steinhaus, compare \cite[\thm 2.5]{RudFA}.
\end{proof}
\begin{lemma}
  \label{lem:liminf}
  $\liminf\limits_{n\to\infty} \lambda_n > 0$.
\end{lemma}
\begin{proof}
  We start by following B\"urger \cite[p.\ 134]{Bue00} and show that
  the spectral radius $\rho(K_\alpha)$ is larger than 1 for
  sufficiently small $\alpha>0$, from which then $\lambda > \alpha >
  0$ follows by \lref{alphalambda}.  Let $J$ be the interval from
  \eref{cusp}.  Then we have
  \begin{displaymath}
    (K_\alpha \ind{J})(x) = 
    \int_J \frac{u(x,y)}{w(y)+\alpha} \dd{y} \ge
    \ind{J}(x) \essinf_{x',y' \in J} u(x',y')
    \int_J (w(y)+\alpha)^{-1} \dd{y}
  \end{displaymath}
  and thus
  \begin{displaymath}
    {\|{K_\alpha}^m\|^{}_1}^{1/m} \ge 
    \essinf_{x,y \in J} u(x,y) \int_J (w(y)+\alpha)^{-1} \dd{y} 
    \qquad\text{for all $m\in\NN$,}
  \end{displaymath}
  which implies for the spectral radius
  \begin{equation}
    \label{eq:rhokalpha}
    \rho(K_\alpha) \ge
    \essinf_{x,y \in J} u(x,y) \int_J (w(y)+\alpha)^{-1} \dd{y} .
  \end{equation}
  The RHS is, as a function of $\alpha$, strictly decreasing.  Thus,
  as a consequence of B.~Levi's monotone convergence theorem
  \cite[\thm III.12.22]{HeSt}, also
  \begin{displaymath}
    \lim_{\alpha\searrow0} \rho(K_\alpha) \ge
    \essinf_{x,y \in J} u(x,y) \int_J (w(y))^{-1} \dd{y} > 1
  \end{displaymath}
  according to \eref{cusp} (including divergence of both sides).
  
  Now we choose $\alpha>0$ such that the RHS of \eref{rhokalpha} is
  larger than or equal to $1+\eps$, with a sufficiently small
  $\eps>0$.  Furthermore, we pick, according to the convergence of the
  quadrature, an $n_0$ with $\essinf_{x,y \in J} u(x,y)
  |Q_n(w+\alpha)^{-1}-Q(w+\alpha)^{-1}| < \eps/2$ for all $n \ge n_0$.
  This way
  \begin{align*}
    (K_{\alpha,n} \ind{J})(x) &= 
    Q_n\bigl(u(x,.) (w+\alpha)^{-1} \ind{J}\bigr) \ge 
    \ind{J}(x) \essinf_{x',y' \in J} u(x',y') Q_n(w+\alpha)^{-1} \\
    &\ge \ind{J}(x) \left(\essinf_{x',y' \in J} u(x',y') Q(w+\alpha)^{-1} -
      \frac{\eps}{2}\right) \ge 
    \ind{J}(x) \left(1+\frac{\eps}{2}\right) \es.
  \end{align*}
  Hence, by \lref{alphalambda}, $\lambda_n > \alpha > 0$ for all $n \ge n_0$,
  from which the claim follows.
\end{proof}

\begin{proof}[Proof of \thref{ieigvalvec}]
  By \lemmas \ref{lem:limsup} and \ref{lem:liminf}, the sequence
  $(\lambda_n)^{}_n$ has a convergent subsequence
  $(\lambda_{n_i})^{}_i$ with limit $\lambda' \in \ocint{0,M}$.
  Consider $(K_{\lambda'} f)(x) = Q(k_{\lambda'}(x,.)f)$ as well as
  $(K_n f)(x) := (K_{\lambda_n,n} f)(x) = Q_n (T+\lambda_n)^{-1}
  (T+\lambda') (k_{\lambda'}(x,.)f)$.  We first show that the
  `distorted' quadrature $\tilde{Q}_{n_i} = Q_{n_i}
  (T+\lambda_{n_i})^{-1} (T+\lambda')$ is convergent.  Note that, for
  $i_0$ large enough, such that $\inf_{j \ge i_0} \lambda_{n_j} > 0$,
  and $i \ge i_0$,
  \begin{equation}
    \label{eq:distnorm}
    \begin{aligned}
      \|(T+\lambda_{n_i})^{-1} (T+\lambda) - 1\|^{}_\infty &=
      \sup_{\|f\|^{}_\infty\le1} 
      \left\|\frac{w+\lambda}{w+\lambda_{n_i}}f - f\right\|^{}_\infty \\
      &\le \|(w + \inf_{j \ge i_0} \lambda_{n_j})^{-1}\|^{}_\infty\, 
      |\lambda-\lambda_{n_i}|\, \|f\|^{}_\infty
      \to 0 \es.
    \end{aligned}
  \end{equation}
  Then, since $(Q_n)$ is convergent by assumption, we have, for all $f
  \in \C(I)$,
  \begin{displaymath}
    \|(\tilde{Q}_{n_i}-Q)f\|^{}_\infty \le
    \|Q_{n_i} ((T+\lambda_{n_i})^{-1} (T+\lambda') - 1) f\|^{}_\infty +
    \|(Q_{n_i}-Q)f\|^{}_\infty \to 0,
  \end{displaymath}
  where the first term vanishes in the limit due to $\sup_m \|Q_m\| <
  \infty$ and \eref{distnorm}.
  
  With this it follows from \thref{nystroem} that $\rho(K_n) = 1$ is also an
  eigenvalue of $K_{\lambda'}$ going with a non-negative eigenfunction.  The
  latter is even a.e.\ positive since, due to the irreducibility
  \eref{irreduc} of $K_{\lambda'}$, there cannot be a set with positive
  measure on which a non-negative eigenfunction vanishes.\footnote{Let
    $\tilde{q}$ be the eigenfunction and $J = \{x : \tilde{q}(x)>0\}$ with $0
    < |J|$.  Assume $|J| < |I|$.  Then, for $x \in I \wo J$, we have $0 =
    \tilde{q}(x) = \int_J k_{\lambda'}(x,y) \tilde{q}(y) \dd{y}$, which
    implies, for a.e.\ $y \in J$, that $k_{\lambda'}(x,y) \tilde{q}(y) = 0$
    and thus $u(x,y) = 0$, contradicting \eref{irreduc}.}  But since,
  according to \thref{existunique}, there is, up to normalization, only one
  positive eigenfunction, we have $\lambda' = \lambda$.  Therefore every
  convergent subsequence of $(\lambda_n)^{}_n$ converges to $\lambda$, and
  thus, due to the boundedness, also the sequence itself.  This proves part
  \ref{it:thmieigval}.

  Along the same line of reasoning, $(a_n^{} q_n^{})$, with $a_n^{} =
  1/\|q_n^{}\|_\infty^{}$, has a convergent subsequence with an a.e.\ positive
  limit function, which equals $aq$ with $a = 1/\|q\|_\infty^{}$.  Therefore,
  $\|a_n^{} q_n^{} - aq\|_\infty^{} \to 0$ for $n\to\infty$.  Now consider
  \begin{align*}
    \|a_n^{} p_n^{} - a p\|^{}_\infty =
    \,&\|a_n^{} p_n^{} - (T+\lambda)^{-1} a q\|^{}_\infty \\
    = \,&\maxp_{k\in\st{N}_n} \sup_{x \in I_{n,k}}
    \bigl|(w(t_{n,k}^{})+\lambda_n)^{-1} a_n^{} q_{n,k}^{} - 
    (w(x)+\lambda)^{-1} a q(x)\bigr| \\
    \le \,&\maxp_{k\in\st{N}_n} \bigl|(w(t_{n,k})+\lambda_n^{})^{-1} - 
    (w(t_{n,k}^{})+\lambda)^{-1}\bigr| \, a_n^{} q_n^{}(t_{n,k}) \, + \\
    &\maxp_{k\in\st{N}_n} (w(t_{n,k}^{})+\lambda)^{-1} 
    \bigl|a_n^{} q_n^{}(t_{n,k}) - a q(t_{n,k}^{})\bigr| \, + \\
    &\maxp_{k\in\st{N}_n} \sup_{x \in I_{n,k}}
    \bigl|(w(t_{n,k}^{})+\lambda)^{-1} a q(t_{n,k}^{}) - 
    (w(x)+\lambda)^{-1} a q(x)\bigr| \es.
  \end{align*}
  The first term is bounded from above by
  \begin{displaymath}
    |\lambda-\lambda_n| \, 
    \|(w + \inf_{m \ge n_0} \lambda_m)^{-1}
    (w+\lambda)^{-1}\|^{}_\infty \es,
  \end{displaymath}
  for $n \ge n_0$ with sufficiently large $n_0$, and vanishes for $n\to\infty$
  due to $\lambda_n \to \lambda$.  The second term vanishes due to the uniform
  convergence of the $a_n^{} q_n^{}$ towards $a q$, and the third due to the
  uniform continuity of $(w+\lambda)^{-1} q$ and \lref{fine}.  With this,
  $a_n^{} p_n^{} \to a p$ in $\LL{\infty}(I)$ and thus in $\LL{1}(I)$.  Hence,
  $a_n^{} \to a$ and $p_n^{} \to p$ in $\LL{1}(I)$, which proves part
  \ref{it:thmieigvec}.
\end{proof}

\section{Discretization---unbounded interval}
\label{sec:unbounded}

Now we assume the types to be taken from $I = \RR$ and the functions
$r$ and $u$ to be continuous.  It will be one aim of this section to
analyze what further conditions have to be imposed in order to allow
for a discretization procedure similar to the one in the previous
section.  In order to do so, we start by a summary of the relevant
theory.

\subsection{The Galerkin method}
\label{sec:galerkin}

In the Galerkin method, an approximation of compact operators is
achieved using projections to finite-dimensional subspaces.  This
method has been reviewed, \eg, by Krasnosel'skii \etal\ \cite[\sect
18]{Kra}.  The results needed in the sequel are collected in
\begin{theorem}
  \label{thm:krasno}
  Let\/ $K$ be a compact linear operator on the Banach space\/ $Y.$
  Consider the eigenvalue equation
  \begin{equation}
    \label{eq:reigvaleqk}
    (K - \nu) g = 0 \es,
  \end{equation}
  which is to be approximated.  To this end, let\/ $(Y_n)$ be a
  sequence of closed subspaces of\/ $Y$ with bounded projections\/
  $P_n$ onto them.  On these subspaces, let the compact linear
  operators\/ $K_n$ be defined, together with the eigenvalue equations
  \begin{equation}
    \label{eq:reigvaleqkn}
    (K_n - \nu_n) g_n = 0 \es.
  \end{equation}
  Assume that
  \begin{equation}
    \label{eq:krasnoassum}
    \|K_n - P_n K\|^{}_{Y_n} \to 0 \es,\quad
    \|K - P_n K\|^{}_Y \to 0 \qquad\text{as $n\to\infty$.}
  \end{equation}
  Then the following statements are true:
  \begin{enumerate}\renewcommand\labelenumi{\textup{(\alph{enumi})}}
  \item\itemlabel{it:eigval} For every\/ $\nu \neq 0$ from
    \eref{reigvaleqk} there is a sequence\/ $(\nu_n)$ of eigenvalues
    of \eref{reigvaleqkn} such that\/ $\nu_n \to \nu$ as\/
    $n\to\infty$.  Conversely, every non-zero limit point of any
    sequence\/ $(\nu_n)$ of eigenvalues of \eref{reigvaleqkn} is an
    eigenvalue of \eref{reigvaleqk}.
  \item\itemlabel{it:eigvec} Every bounded sequence\/ $(g_n)$ of
    eigenvectors of \eref{reigvaleqkn} associated with eigenvalues\/
    $\nu_n \to \nu \neq 0$ contains a convergent subsequence; the
    limit of any convergent subsequence\/ $(g_{n_i})^{}_i$ is an
    eigenvector of \eref{reigvaleqk} associated with the eigenvalue\/
    $\nu$ \pr{unless the limit is zero}.
  \end{enumerate}
\end{theorem}
\begin{proof}
  See \cite[\thms 18.1, 18.2]{Kra}.
\end{proof}

A sufficient condition for the validity of the second assumption in
\eref{krasnoassum} is given by
\begin{proposition}
  \label{prop:pointwise}
  Let\/ $X$ be a normed space, $Y$ a Banach space, and\/ $K\colon X \to Y$ a
  compact linear operator.  For bounded linear operators\/ $P_n \colon Y \to
  Y$ $(n\in\NN)$ with\/ $P_n \to \Id$ pointwise for\/ $n\to\infty$, the
  operators\/ $P_n K$ approximate\/ $K$, \ie $\|P_n K - K\| \to 0$.
\end{proposition}
\begin{proof}
  Follow the proof of \cite[\satz II.3.5]{Wer}, where the additional
  assumptions on $X$ and $(P_n)$ are not used.
\end{proof}

\subsection{Application to kernel operators}

In our case of the COA model we have $X = Y = \LL{1}(\RR)$ and $K$ is
of the form
\begin{equation}
  \label{eq:kernelopr}
  (Kf)(x) = \int_\RR k(x,y) f(y) \dd{y}
  \qquad\text{for all $x \in \RR$}
\end{equation}
with a measurable kernel $k \colon \RR\times\RR \to \RR$.  Therefore, for the
Galerkin method to work, it is necessary that, for $\LL{1}(\RR)$, operators
$P_n$ as in \pref{pointwise} exist.  We will explicitly construct such
operators using a sequence $(\{I_{n,k} : 1 \le k \le N_n\})^{}_n$ of families
of disjoint intervals that get finer and finer and also ultimately cover every
bounded interval.\footnote{Both properties are formally captured by
  \ref{it:cover} in \pref{apprprop}.}
\begin{proposition}
  \label{prop:apprprop}
  Let\/ $Y$ be the Banach space\/ $\LL{1}(\RR)$ and finite-dimensional
  subspaces\/ $Y_n$ of\/ $Y$ chosen to consist of all step functions with
  prescribed \pr{bounded} intervals\/ $I_{n,k}$ \pr{$k \in \st{N}_n :=
    \{1,\ldots,N_n\}$} with the following properties:
  \begin{itemize}  
  \item[\itlabel{(I1)}{it:cover}] For every bounded interval\/
    $I\subset\RR$ and every\/ $\eps>0$ there is an\/ $n_0$ such that,
    for all\/ $n \ge n_0$, a set\/ $L\subset\st{N}_n$ exists for
    which\/ $I_{n,L} := \bigcup_{\ell \in L} I_{n,\ell}$ satisfies\/
    $|I \wo I_{n,L}| = 0$ and\/ $|I_{n,L} \wo I| < \eps$.  \pr{We then
    say that\/ $I$ is \textdef{$\eps$-optimally covered}.}
  \item[\itlabel{(I2)}{it:disjoint}] $|I_{n,k} \cap I_{n,\ell}| = 0$
    for all\/ $n\in\NN$ and\/ $1 \le k < \ell \le N_n$.
  \end{itemize}
  Then, with the characteristic functions\/ $\varphi_{n,k} =
  \ind{I_{n,k}}$, the projections\/ $P_n$ onto the subspaces\/ $Y_n$
  spanned by\/ $\{\varphi_{n,k} : k \in \st{N}_n\}$ are given by
  \begin{displaymath}
    P_n f = \sum_{k=1}^{N_n} \varphi_{n,k} \frac{1}{|I_{n,k}|} 
      \int_{I_{n,k}} \! f(x) \dd{x}
    \qquad\text{for\/ $f\in\LL{1}(\RR)$,}
  \end{displaymath}
  where\/ $\int_{I_{n,k}} f(x) \dd{x}$ are \textdef{conditional expectations}
  \pr{compare \textup{\cite[\thm IV.2.4]{Schae}}\footnote{See also \cite{Diss}
      for a discussion of the connection to the approximation property of
      Banach spaces.}}.  The projections satisfy\/ $\|P_n\|^{}_1 = 1$ and\/
  $P_n \to \Id$ pointwise.
\end{proposition}
\begin{proof}
  Obviously, the subspaces $Y_n$ are closed, finite-dimensional, and
  the $P_n$ are, due to \ref{it:disjoint}, projections onto them.
  Since
  \begin{displaymath}
    \|P_n f\| = \sum_{k=1}^{N_n} \int_{I_{n,k}} |f(x)|\dd{x} \le 
    \int_\RR |f(x)|\dd{x} = \|f\|^{}_1 
    \qquad\text{for every $f \in \LL{1}(\RR)$}
  \end{displaymath}
  and $\|P_n \varphi_{n,k}\| = \|\varphi_{n,k}\|$ for every $k \in
  \st{N}_n$, we have $\|P_n\|=1$.
  
  We now show that $P_n \to \Id$ pointwise.  To this end, let $f \in
  \LL{1}(\RR)$ and $\eps>0$ be given.  Remember that the set of all
  step functions is, by definition, dense in $\LL{1}(\RR)$, compare
  \cite[\sect VI.3]{Lan}.  Therefore, we can find a step function
  $\psi = \sum_{k=1}^m \psi_k \ind{J_k}$ (with bounded intervals
  $J_k$) that satisfies $\|f-\psi\|^{}_1 < \eps/3$.  Due to
  \ref{it:cover} we can now choose an $n_0$ such that
  $|\bigcup_{k=1}^m J_k \wo \bigcup_{k=1}^{N_n} I_{n,k}| = 0$ for all
  $n \ge n_0$.  Then, the only contributions to $\|P_n \psi -
  \psi\|^{}_1$ are due to mismatches at the boundaries of the $J_k$.
  Therefore, let $J^+_k$ and $J^-_k$ ($k \in \st{N}_n$) be open
  intervals of measure $\eta = \eps / (12m \max_k |\psi_k|)$ that
  contain the right and left boundary points of $J_k$, respectively.
  Choosing $n_1 \ge n_0$ according to \ref{it:cover} large enough such
  that every $J^\pm_k$ is $\eta$-optimally covered for $n \ge n_1$, we
  have $\|P_n \psi - \psi\|^{}_1 < 2 \sum_{k=1}^m 2\eta \psi_k \le
  \eps/3$ for $n \ge n_1$.  Putting everything together yields, for $n
  \ge n_1$,
  \begin{displaymath}
    \|P_n f - f\|^{}_1 \le 
    \|P_n (f-\psi)\|^{}_1 + \|P_n \psi - \psi\|^{}_1 +
      \|\psi - f\|^{}_1 < \eps \es,
  \end{displaymath}
  which proves $\|P_n f - f\| \to 0$ for $n\to\infty$ and thus the
  approximation property.
\end{proof}

With respect to a kernel operator $K$ of the form \eref{kernelopr} and
some $f = \sum_{k=1}^{N_n} \varphi_{n,k} f_k$ in $Y_n$, the above
procedure amounts to the discretization
\begin{align*}
  (P_n K f) &= 
  \sum_{k=1}^{N_n} \varphi_{n,k} \frac{1}{|I_{n,k}|} 
    \left( \int_{I_{n,k}}
    \int_\RR k(x,y) \sum_{\ell=1}^{N_n} f_\ell \varphi_{n,\ell}(y) 
    \dd{y}\dd{x} \right) \\
  &= \sum_{k=1}^{N_n} \varphi_{n,k} \sum_{\ell=1}^{N_n} \frac{1}{|I_{n,k}|}
    \left( \int_{I_{n,k}} \int_{I_{n,\ell}} k(x,y) \dd{y}\dd{x} 
  \right) f_\ell =:
  \sum_{k=1}^{N_n} \varphi_{n,k} \sum_{\ell=1}^{N_n} M_{n,k\ell} f_\ell
\end{align*}
with an $N_n \times N_n$ matrix $\mt{M}_n = (M_{n,k\ell})$.  The
corresponding eigenvalue equation is
\begin{displaymath}
  \mt{M}_n^{} \vc{g}_n^{} = \nu_n^{} \vc{g}_n^{} \es,
  \qquad\text{or equivalently}\quad 
  P_n^{} K g_n^{} = \nu_n^{} g_n^{} \es,
\end{displaymath}
where $g_n^{} \in Y_n^{}$ is granted due to the projection property.
An example of intervals $I_{n,k}$ satisfying \ref{it:cover} and
\ref{it:disjoint} is $I_{n,k} = \ccint{-n+2^{-n}(k-1),-n+2^{-n}k}$
with $k \in \st{N}_n = \{1,\ldots,2^{n+1}n\}$.

With respect to compactness of $K$, following J\"orgens \cite[\sects 11,
12]{Joe}, we extract
\begin{proposition}
  \label{prop:kcomp}
  A kernel operator\/ $K$ on\/ $\LL{1}(\RR)$ of the form \eref{kernelopr} is
  compact if it satisfies the following conditions:
  \begin{itemize}
  \item[\itlabel{(C1)}{it:comp1}] The function\/ $x \to k(x,.)$ from\/
    $\RR$ to $\LL{1}(\RR)$ is continuous and bounded.
  \item[\itlabel{(C2)}{it:comp2}] For every\/ $\eps>0$ there exists a
    finite open covering\/ $(V_1,\ldots,V_n)$ of\/ $\RR$ and points\/
    $x_j \in V_j$ such that\/ $\|k(x,.)-k(x_j,.)\|^{}_1 < \eps$ for
    all\/ $x \in V_j$ and all\/ $j$.
  \item[\itlabel{(C3)}{it:comp3}] The function\/ $y \to k(.,y)$ from\/
    $\RR$ to $\LL{1}(\RR)$ is continuous and bounded.
  \item[\itlabel{(C4)}{it:comp4}] For every\/ $\eps>0$ there exists a
    finite open covering\/ $(W_1,\ldots,W_n)$ of\/ $\RR$ and points\/
    $y_j \in W_j$ such that\/ $\|k(.,y)-k(.,y_j)\|^{}_1 < \eps$ for
    all\/ $y \in W_j$ and all\/ $j$.
  \end{itemize}
\end{proposition}
\begin{proof}
  First, as in \cite[\sect 12.4]{Joe}, we consider the dual system
  $\scal{\C(\RR), \CI(\RR)}$ with the bilinear form $\scal{f,g} = \int_\RR
  f(x) g(x) \dd{x}$.  Here, $\C(\RR)$ is equipped with the supremum norm
  $\|.\|^{}_\infty$ and $\CI(\RR) := \C(\RR) \cap \LL{1}(\RR)$ with the norm
  $\tnorm{.} := \max\{\|.\|^{}_\infty, \|.\|^{}_1\}$.  With this, we define
  the transposed $K\trp$ of $K$ via $(K\trp g)(y) = \int_\RR g(x) k(x,y)
  \dd{x}$, for all $y \in \RR$.  Then, by \ref{it:comp1}--\ref{it:comp4} and
  \cite[\saetze 12.2, 12.3]{Joe}, the compactness of $K$ and $K\trp$ on
  $\C(\RR)$ follows.
  
  As both $K\trp$ and $K$ are bounded as operators on $\C(\RR)$, they are, at
  the same time, Hille--Tamarkin operators in $\HT{\infty\infty}(\RR)$ since
  the respective norm, $\thnorm{.}^{}_{\infty\infty}$ in \eref{hiltampq}, is
  just given by $\sup_{x\in\RR} \int_\RR k(x,y) \dd{y}$ and $\sup_{y\in\RR}
  \int_\RR k(x,y) \dd{x}$, respectively \cite[\saetze 12.2, 12.3]{Joe}.  Then,
  according to \cite[\satz 11.5]{Joe}, $K$ and $K\trp$ can also be regarded as
  bounded operators on $\LL{1}(\RR)$; thus, both map $\CI(\RR)$ into itself.
  Due to \cite[\satz 12.6]{Joe} there is, for every $\eps>0$, an operator of
  finite rank, $K_\eps$, with $\tnorm{K_\eps - K} < \eps$, where $\tnorm{A} :=
  \max\{\|A\|^{}_\infty, \|A\trp\|^{}_\infty\}$ is a norm for the Banach
  algebra of all operators on $\C(\RR)$ that map $\CI(\RR)$ into itself and
  have a transposed of the same kind.  We have $\tnorm{Af} \le \tnorm{A}
  \tnorm{f}$ for $f \in \CI(\RR)$, see \cite[\sect 12.4]{Joe}.  Thus,
  $\tnorm{A}$ can serve as an upper bound for the operator norm of $A$ on
  $\CI(\RR)$.  Therefore, $K$ is compact as an operator on $\CI(\RR)$ and can
  be approximated by $K_\eps$.  Furthermore, according to \cite[\satz
  11.5]{Joe}, $\|K_\eps - K\|^{}_1 \le \tnorm{(K_\eps - K)\trp} < \eps$ holds.
  Hence, $K$ is compact as an operator on $\LL{1}(\RR)$ as well.
\end{proof}

\subsection{Application to the COA model}

Checking the compactness of $K_\alpha$ by conditions
\ref{it:comp1}--\ref{it:comp4} of \pref{kcomp}, we would be able to apply
\thref{krasno} and approximate $K_\alpha$ by operators of finite rank.
However, the original system is described by the (non-compact) operator $A = T
- U$, not by some $K_\alpha$.  It will be shown that it is indeed possible to
discretize the operators $T$ and $U$ directly by applying the projections
$P_n$ from \pref{apprprop}, if further restrictions apply.  Then, even more
generally, the approximation can be done by choosing arbitrary points in the
intervals $I_{n,k}$ at which the functions $w$ and $u$ are sampled.  Both
procedures will now be described in detail.

In the first setting, $K_\lambda$ is approximated by $K_n := P_n U
(P_n T + \lambda_n)^{-1}$.  Explicitly, for $f \in Y_n$ with $f =
\sum_{k=1}^{N_n} f_k \varphi_{n,k}$, it reads
\pagebreak[2]
\begin{alignat*}{2}
  P_n T f &= 
  \mskip5mu\sum_{k=1}^{N_n}\mskip5mu \varphi_{n,k}  
    \frac{1}{|I_{n,k}|} \int_{I_{n,k}} w(x) \dd{x} \, f_k &&= 
  \mskip5mu\sum_{k=1}^{N_n}\mskip5mu \varphi_{n,k} w(t^w_{n,k}) f_k \es,\\
  P_n U f &=
  \sum_{k,\ell=1}^{N_n} \varphi_{n,k}
    \frac{1}{|I_{n,k}|} \int_{I_{n,k}} \int_{I_{n,\ell}}
      u(x,y) \dd{y} \dd{x} \, f_\ell &&=
  \sum_{k,\ell=1}^{N_n} \varphi_{n,k}
    |I_{n,\ell}| u(t^{ux}_{n,k\ell}, t^{uy}_{n,k\ell}) f_\ell \es,
\end{alignat*}
with appropriate points $t^w_{n,k}$, $t^{ux}_{n,k\ell} \in I^{}_{n,k}$
and $t^{uy}_{n,k\ell} \in I^{}_{n,\ell}$ that satisfy
\begin{equation}
  \label{eq:tuxy}
  \frac{1}{|I_{n,k}|} \int_{I_{n,k}} u(x,t^{uy}_{n,k\ell}) \dd{x} =
  u(t^{ux}_{n,k\ell}, t^{uy}_{n,k\ell}) \es.
\end{equation}
These exist due to the continuity of $w$ and $u$.  But more generally,
we may pick the points arbitrarily from the respective intervals.

In either case, we define the $N_n \times N_n$ matrices $\mt{T}_n$,
$\mt{U}_n$, and $\mt{A}_n := \mt{T}_n - \mt{U}_n$ via
\begin{equation}
  \label{eq:rdeftnun}
  T_{n,kk} := w(t^w_{n,k}) \es,\qquad
  U_{n,k\ell} := |I_{n,\ell}| \, u(t^{ux}_{n,k\ell}, t^{uy}_{n,k\ell}) \es.
\end{equation}
The corresponding operators in $Y_n$ are given by
\begin{displaymath}
  T_n f = \sum_{k=1}^{N_n} \varphi_{n,k} T_{n,kk} f_k \es, \qquad
  U_n f = \sum_{k,\ell=1}^{N_n} \varphi_{n,k} U_{n,k\ell} f_\ell \es, \qquad
  A_n = T_n - U_n \es,
\end{displaymath}
again with $f = \sum_{k=1}^{N_n} f_k \varphi_{n,k}$.  For notational
convenience, we also define the matrices $\mt{P}_{\alpha,n}$ by
\begin{equation}
  \label{eq:rdefkan}
  P_n K_\alpha f =
    \sum_{k,\ell=1}^{N_n} \varphi_{n,k}
    \frac{1}{|I_{n,k}|} \int_{I_{n,k}} \int_{I_{n,\ell}}
      \frac{u(x,y)}{w(y)+\alpha} \dd{y} \dd{x} f_\ell =:
    \sum_{k,\ell=1}^{N_n} \varphi_{n,k} P_{\alpha,n,k\ell} f_\ell \es.
\end{equation}
The eigenvalue equation to be solved is
\begin{displaymath}
  (\mt{A}_n^{} + \lambda_n^{}) \vc{p}_n^{} = 0 \es,
\end{displaymath}
which is equivalent to
\begin{equation}
  \label{eq:rdiscr}
  (A_n^{} + \lambda_n^{}) p_n^{} = 0 \es,
\end{equation}
where $p_n^{} = \sum_{k=1}^{N_n} p_{n,k}^{} \varphi_{n,k}^{} \in
Y_n^{}$.  With $K_{\alpha,n} = U_n (T_n + \alpha)^{-1}$, $\alpha >
-\min_{k\in\st{N}_n} w(t_{n,k})$, and $q_n^{} = (T_n^{} +
\lambda_n^{}) p_n^{}$ also the eigenvalue equation
\begin{displaymath}
  (K_{\lambda_n,n} - 1) q_n^{} = 0
\end{displaymath}
is solved (and vice versa), \cf \lref{connectka}.  (The inequality
$\lambda_n > -\min_{k\in\st{N}_n} w(t_{n,k})$ follows from
\thref{existunique}.)

For these procedures to be valid approximations, the first condition
in \eref{krasnoassum}, that is, $\|K_n - P_n K\|^{}_{Y_n} \to 0$, has
to be true for $K = K_\lambda$ and $K_n = U_n (T_n + \lambda_n)^{-1}$.
This, however, is not given automatically.  Problems arise from the
fact that in $K_n$ the averaging defined by $P_n$ (or, more generally,
the sampling) is applied to the enumerator and denominator of
$k_{\lambda_n}$ separately, whereas in $P_n K$ the quotient
$k_\lambda$ is averaged as such.  It turns out that some additional
requirements of uniform continuity are sufficient for the convergence.
This is made precise in the following two propositions.
\begin{proposition}
  \label{prop:fixeda}
  Suppose that the following conditions are true:
  \begin{itemize}
  \item[\itlabel{(S1)}{it:uxunicont}] $u(x,.)$ is uniformly continuous
    for all\/ $x \in \RR$.
  \item[\itlabel{(S2)}{it:kaunicont}] $k_\alpha$ is uniformly continuous on\/
    $I\times\RR$ for all\/ $\alpha>0$ and all bounded\/ $I\subset\RR$.
  \item[\itlabel{(S3)}{it:wmin}] There is a function\/ $\wmin
    \colon \RR \to \RRnn$, satisfying
    \begin{displaymath}
      \int_\RR \sup_{y\in\RR} \frac{u(x,y)}{\wmin(y)+\alpha} \dd{x} < 
      \infty \qquad\text{for all\/ $\alpha>0$,}
    \end{displaymath}
    and an\/ $n_0\in\NN$ such that\/ $w(y) \ge \wmin(y')$ for all\/ $n
    \ge n_0$, $\ell\in\st{N}_n$, and\/ $y$, $y' \in I_{n,\ell}$.
  \end{itemize}
  Then, for\/ $K = K_\alpha$ and\/ $K_n = P_n U (P_n T + \alpha)^{-1}$
  with any\/ $\alpha>0$ and the projections\/ $P_n$ from
  \pref{apprprop}, the first condition in \eref{krasnoassum} is
  fulfilled, \ie $\|K_n - P_n K\|^{}_{Y_n} \to 0$ as $n\to\infty$.
  The same is true for\/ $K_n = K_{\alpha,n} = U_n (T_n +
  \alpha)^{-1}$ with the more general discretization from above if in
  addition to \ref{it:uxunicont}--\ref{it:wmin} the following
  condition is satisfied:
  \begin{itemize}
  \item[\itlabel{(S4)}{it:umax}] There is a function\/ $\umax
    \colon \RR\times\RR \to \RRnn$, satisfying
    \begin{displaymath}
      \int_\RR \sup_{y\in\RR} \frac{\umax(x,y)}{\wmin(y)+\alpha} \dd{x} < 
      \infty \qquad\text{for all\/ $\alpha>0$,}
    \end{displaymath}
    and an\/ $n_1 \ge n_0$ such that\/ $u(x,y) \le \umax(x',y)$ for
    all\/ $n \ge n_1$, $k\in\st{N}_n$, $y\in\RR$, and\/ $x$, $x' \in
    I_{n,k}$.
  \end{itemize}
\end{proposition}
Let us split the rather technical proof into a couple of digestible
lemmas.
\begin{lemma}
  \label{lem:inner}
  If conditions \ref{it:uxunicont} and \ref{it:kaunicont} are true, then
  for every\/ $\eps>0$ and every compact interval\/ $I \subset \RR$
  there is an\/ $n_2$ such that for all\/ $n \ge n_2$ and all $k$,
  $\ell \in \st{N}_n$ with\/ $I_{n,k} \cap I \neq \emptyset$ we have
  \begin{displaymath}
    \frac{1}{|I_{n,\ell}|}
    \left| P_{\alpha,n,k\ell} - 
      \frac{U_{n,k\ell}}{T_{n,\ell\ell}+\alpha} \right|
    < \frac{\eps}{|I|} \es.
  \end{displaymath}
\end{lemma}
\begin{proof}
  Let $\eps$ and $I$ be given as above and $I_0 = \bigcup_{n\in\NN} \bigcup_{k
    : I_{n,k} \cap I \neq \emptyset} I_{n,k}$, which is a bounded interval due
  to \ref{it:cover} from \pref{apprprop}.  By assumptions \ref{it:uxunicont}
  and \ref{it:kaunicont}, $u$ and $k_\alpha$ are uniformly continuous on
  $\overline{I_0}\times\RR$.  Further, $(w+\alpha)^{-1}$ is bounded by
  $\alpha^{-1}$.  Thus, there is an $n_2$ such that, for every $n \ge n_2$ and
  $k$, $\ell\in\st{N}_n$ with $I_{n,k} \cap I \neq \emptyset$,
  \begin{multline*}
    \left| \frac{1}{|I_{n,k}|} \frac{1}{|I_{n,\ell}|} 
      \int_{I_{n,k}} \int_{I_{n,\ell}} 
      \frac{u(x,y)}{w(y)+\alpha} \dd{y} \dd{x} -
      \frac{u(t^{ux}_{n,k\ell}, t^{uy}_{n,k\ell})}{w(t^w_{n,\ell})+\alpha}
    \right| \\
    = \left|
      \frac{u(t^{kx}_{n,k\ell}, t^{ky}_{n,k\ell})}{w(t^{ky}_{n,k\ell})+\alpha}-
      \frac{u(t^{ux}_{n,k\ell}, t^{uy}_{n,k\ell})}{w(t^w_{n,\ell})+\alpha}
    \right| \\
    \le \left|
      \frac{u(t^{kx}_{n,k\ell}, t^{ky}_{n,k\ell})}{w(t^{ky}_{n,k\ell})+\alpha}-
      \frac{u(t^{ux}_{n,k\ell}, t^w_{n,\ell})}{w(t^w_{n,\ell})+\alpha}
    \right| +
    \frac{|u(t^{ux}_{n,k\ell},t^w_{n,\ell}) - 
        u(t^{ux}_{n,k\ell},t^{uy}_{n,k\ell})|}
      {w(t^w_{n,\ell})+\alpha} <
    \frac{\eps}{|I|} \es.
  \end{multline*}
  Here, the points $t^{kx}_{n,k\ell} \in I_{n,k}$ and
  $t^{ky}_{n,k\ell} \in I_{n,\ell}$ are chosen such that the first
  equality holds, which is possible due to the continuity of
  $k_\alpha$.  From this the claim follows easily with \eref{rdeftnun}
  and \eref{rdefkan}.
\end{proof}
\begin{lemma}
  \label{lem:outer1}
  For every\/ $\eps>0$ there is a compact interval\/ $I_1$ such that,
  for all intervals\/ $I \supset I_1$ and all\/ $n\in\NN$,
  \begin{displaymath}
    \sum_{\substack{k \\ I_{n,k} \cap I = \emptyset}} |I_{n,k}|
      \max_{\ell\in\st{N}_n} \frac{P_{\alpha,n,k\ell}}{|I_{n,\ell}|} < 
    \eps \es.
  \end{displaymath}
\end{lemma}
\begin{proof}
  Due to \ref{it:u4} there is a compact interval $I_1$ such that, for
  all $I \supset I_1$,
  \begin{align*}
    \sum_{\substack{k \\ I_{n,k} \cap I = \emptyset}} |I_{n,k}|
      \max_{\ell\in\st{N}_n} \frac{P_{\alpha,n,k\ell}}{|I_{n,\ell}|} \le
    \sum_{\substack{k \\ I_{n,k} \cap I = \emptyset}} |I_{n,k}|
      \frac{1}{|I_{n,k}|} \int_{I_{n,k}} \max_{y\in\RR}
        \frac{u(x,y)}{w(y)+\alpha} \dd{x} \\
    \le \int_{\RR \wo I_1} \max_{y\in\RR} \frac{u(x,y)}{w(y)+\alpha} \dd{x} < 
    \eps \es,
  \end{align*}
  which proves the claim.
\end{proof}
\begin{lemma}
  \label{lem:outer2}
  If condition \ref{it:wmin} is true, and if
  \begin{enumerate}\renewcommand\labelenumi{\textup{(\roman{enumi})}}
  \item\itemlabel{it:outer2a} $U_{n,k\ell} =
    \frac{|I_{n,\ell}|}{|I_{n,k}|} \int_{I_{n,k}}
    u(x,t^{uy}_{n,k\ell}) \dd{x}$ for all\/ $k$, $\ell\in\st{N}_n$ or
  \item\itemlabel{it:outer2b} condition \ref{it:umax} is fulfilled,
  \end{enumerate}
  then for every\/ $\eps>0$ there is a compact interval\/ $I_2$ such
  that, for all intervals\/ $I \supset I_2$ and all\/ $n\in\NN$,
  \begin{displaymath}
    \sum_{\substack{k \\ I_{n,k} \cap I = \emptyset}} |I_{n,k}|
      \max_{\ell\in\st{N}_n} 
      \frac{U_{n,k\ell}}{|I_{n,\ell}| (T_{n,\ell\ell}+\alpha)} < \eps \es.
  \end{displaymath}
\end{lemma}
\begin{proof}
  In case \ref{it:outer2a} we have, using \eref{tuxy},
  \begin{multline*}
    \sum_{\substack{k \\ I_{n,k} \cap I = \emptyset}} |I_{n,k}|
      \max_{\ell\in\st{N}_n}
      \frac{u(t^{ux}_{n,k\ell}, t^{uy}_{n,k\ell})}{w(t^w_{n,\ell})+\alpha} =
    \sum_{\substack{k \\ I_{n,k} \cap I = \emptyset}}
      \max_{\ell\in\st{N}_n} 
      \frac{\int_{I_{n,k}} u(x, t^{uy}_{n,k\ell}) \dd{x}}
        {w(t^w_{n,\ell})+\alpha} \\
    \le \sum_{\substack{k \\ I_{n,k} \cap I = \emptyset}} \int_{I_{n,k}} 
      \max_{y\in\RR} \frac{u(x,y)}{\wmin(y)+\alpha} \dd{x} \le 
    \int_{\RR \wo I_2} \max_{y\in\RR} 
      \frac{u(x,y)}{\wmin(y)+\alpha} \dd{x} < \eps
  \end{multline*}
  for some compact interval $I_2$, due to \ref{it:wmin}, and 
  all intervals $I \supset I_2$.  In case \ref{it:outer2b} we can
  find, due to \ref{it:umax}, a compact interval $I_2$ such that, for
  all intervals $I \supset I_2$,
  \begin{multline*}
    \sum_{\substack{k \\ I_{n,k} \cap I = \emptyset}} |I_{n,k}|
      \max_{\ell\in\st{N}_n}
        \frac{u(t^{ux}_{n,k\ell}, t^{uy}_{n,k\ell})}{w(t^w_{n,\ell})+\alpha}\le
    \sum_{\substack{k \\ I_{n,k} \cap I = \emptyset}} |I_{n,k}|
      \max_{y\in\RR} \frac{u(t^{ux}_{n,k\ell}, y)}{\wmin(y)+\alpha} \\
    \le \sum_{\substack{k \\ I_{n,k} \cap I = \emptyset}} \max_{y\in\RR} 
      \frac{\int_{I_{n,k}} \umax(x,y) \dd{x}}{\wmin(y)+\alpha} \le
    \int_{\RR \wo I_2} \max_{y\in\RR} 
      \frac{\umax(x,y)}{\wmin(y)+\alpha} \dd{x} < \eps \es.
  \end{multline*}
  \Tctp
\end{proof}
\begin{proof}[Proof of \pref{fixeda}]
  Let $\eps>0$ be given.  Choose a compact interval $I$ such that $I
  \supset I_1 \cup I_2$ with $I_1$ and $I_2$ from \lemmas
  \ref{lem:outer1} and \ref{lem:outer2}.  Let $I_3 =
  \overline{\bigcup_{\substack{n,k : I_{n,k} \cap I \neq \emptyset}}
  I_{n,k}}$.  Further, let $n_0$ be as in \ref{it:wmin}, $n_1$ as in
  \ref{it:umax} (or $n_0 = n_1$ if not applicable), $n_2$ as in
  \lref{inner}, and $n \ge \max\{n_0, n_1, n_2\}$.  Then
  \begin{multline*}
    \| P_n K_\alpha - K_{\alpha,n} \|^{}_{Y_n} = 
    \sup_{\substack{f \in Y_n \\ \|f\|^{}_{Y_n}\le1}} 
       \sum_{k=1}^{N_n} |I_{n,k}|
    \left| \sum_{\ell=1}^{N_n} \left( 
        P_{\alpha,n,k\ell} - \frac{U_{n,k\ell}}{T_{n,\ell\ell}+\alpha} \right)
      f_\ell \right| \\
    \le \biggl( \sum_{\substack{k \\ I_{n,k} \cap I \neq \emptyset}} +
      \sum_{\substack{k \\ I_{n,k} \cap I = \emptyset}} \biggr)
    |I_{n,k}| \max_{\ell\in\st{N}_n} \frac{1}{|I_{n,\ell}|}
    \left| P_{\alpha,n,k\ell} - 
      \frac{U_{n,k\ell}}{T_{n,\ell\ell}+\alpha} \right| \\
    \le \sum_{\substack{k \\ I_{n,k} \cap I \neq \emptyset}} 
      |I_{n,k}| \frac{\eps}{|I_3|} +
    \sum_{\substack{k \\ I_{n,k} \cap I = \emptyset}}
    |I_{n,k}| \left( \max_{\ell\in\st{N}_n} 
        \frac{P_{\alpha,n,k\ell}}{|I_{n,\ell}|} +
      \max_{\ell\in\st{N}_n} 
        \frac{U_{n,k\ell}}{|I_{n,\ell}| (T_{n,\ell\ell}+\alpha)}
    \right)
    < 3\eps
  \end{multline*}
  according to \lemmas \ref{lem:inner}--\ref{lem:outer2}.  From this
  the claim follows.
\end{proof}
\begin{proposition}
  \label{prop:convergenta}
  Let\/ $\alpha_n > -\min_{k\in\st{N}_n} w(t_{n,k})$ with\/ $\alpha_n
  \to \alpha > 0$ as\/ $n\to\infty$ and the hypotheses of
  \pref{fixeda} be satisfied.  Then\/ $\| K_{\alpha_n,n} - P_n
  K_\alpha \|^{}_{Y_n} \to 0$.
\end{proposition}
\begin{proof}
  Consider
  \begin{multline*}
    \| P_n K_\alpha - U_n (T_n + \alpha_n)^{-1} \|^{}_{Y_n} \\
    \le \| P_n K_\alpha - U_n (T_n + \alpha)^{-1} \|^{}_{Y_n} +
       \| U_n [ (T_n + \alpha_n)^{-1} -
         (T_n + \alpha)^{-1} ] \|^{}_{Y_n} \es.
  \end{multline*}
  The first term tends to zero as $n\to\infty$ according to
  \pref{fixeda}.  For the second, choose $n_0$ such that $\inf_{n \ge
    n_0} \alpha_n > 0$.  Then, for $n \ge n_0$,
  \begin{multline*}
    \| U_n [ (T_n + \alpha_n)^{-1} - 
      (T_n + \alpha)^{-1} ] \|^{}_{Y_n} \\
    = |\alpha-\alpha_n| \,
      \| U_n (T_n + \alpha_n)^{-1} 
        (T_n + \alpha)^{-1} \|^{}_{Y_n} \le
    |\alpha-\alpha_n| \, \|U\|^{}_Y (\inf_{n \ge n_0} \alpha_n)^{-1} 
      \alpha^{-1} \es.
  \end{multline*}
  This vanishes as $n\to\infty$ since all constants that occur are
  finite, from which the claim follows.
\end{proof}

\subsection{Convergence of eigenvalues and eigenvectors}

Let us now show
\begin{theorem}
  \label{thm:reigvalvec}
  With the notation and assumptions from \sref{operatornot} and\/
  $\lambda$, $p$, $\lambda_n$, $p_n^{}$ as in \eref{coaequilop} and
  \eref{rdiscr}, we have
  \begin{enumerate}\renewcommand\labelenumi{\textup{(\alph{enumi})}}
  \item\itemlabel{it:thmreigval} $\lim_{n\to\infty} \lambda_n =
    \lambda > 0$ and
  \item\itemlabel{it:thmreigvec} $\lim_{n\to\infty} \|p_n^{}-p\|^{}_1
    = 0$, \ie the probability measures corresponding to these
    densities converge in total variation.
  \end{enumerate}
\end{theorem}
The plan is the same as described in \sref{ieigvalvec}.  The proofs,
however, are quite different due to the more general setup.
\begin{lemma}
  \label{lem:rlimsup}
  There is a constant\/ $M>0$ such that\/ $\limsup_{n\to\infty}
  \lambda_n \le M$.
\end{lemma}
\begin{proof}
  Choose an $\alpha>0$ such that $\|K_\alpha\|^{}_Y \le 1-\eps$ for
  some $0 < \eps < 1$, which is possible since $\|K_\alpha\|^{}_Y \to
  0$ for $\alpha \to \infty$.  Then, for all $n \ge n_0$ with some
  $n_0$, due to \propositions \ref{prop:pointwise} and
  \ref{prop:fixeda}, $| \|P_n K_\alpha\|^{}_Y - \|K_\alpha\|^{}_Y |\le
  \eps/3$ and $| \|K_{\alpha,n}\|^{}_{Y_n} - \|P_n K_\alpha\|^{}_{Y_n}
  | \le \eps/3$.  For these $n$, we have
  \begin{align*}
    \rho(K_{\alpha,n}) \le \|K_{\alpha,n}\|^{}_{Y_n} 
    \le \|P_n K_\alpha\|^{}_{Y_n} + \eps/3
    &\le \|P_n K_\alpha\|^{}_Y + \eps/3 \\
    &\le \|K_\alpha\|^{}_Y + 2\eps/3 
    \le 1 - \eps/3 < 1
  \end{align*}
  and thus $\lambda_n < \alpha$ by \lref{alphalambda}.  Then, with $M
  = \alpha$, the claim follows.
\end{proof}
\begin{lemma}
  \label{lem:rliminf}
  $\liminf_{n\to\infty} \lambda_n > 0$.
\end{lemma}
\begin{proof}
  In a modification of the proof of \lref{liminf}, we choose
  $\alpha>0$ such that $\rho(K_\alpha) \ge 1+\eps$ with a sufficiently
  small $\eps>0$.  We know from the theorem of Jentzsch \cite[\thm
  V.6.6]{Schae} that $\rho(K_\alpha)$ is a simple eigenvalue of
  $K_\alpha$ and the only one with a positive eigenfunction.  The same
  is true for $\rho(K_{\alpha,n})$ \wrt $K_{\alpha,n}$ (as an operator
  in $Y_n$).  \Thref{krasno} together with \pref{fixeda} implies that
  there is a sequence of eigenvalues $\nu_n$ of $K_{\alpha,n}$ with
  limit $\rho(K_\alpha)$.  Therefore, $\liminf_{n\to\infty}
  \rho(K_{\alpha,n}) \ge \rho(K_\alpha) \ge 1+\eps$ and thus
  $\lambda_n > \alpha > 0$ for sufficiently large $n$.  From this the
  claim follows.
\end{proof}
\begin{proof}[Proof of \thref{reigvalvec}]
  From \lref{rlimsup} and \ref{lem:rliminf} we conclude that there is a
  convergent subsequence $(\lambda_{n_i})^{}_i$ with limit $\lambda' \in
  \ocint{0,M}$.  Then, due to \pref{convergenta}, $K_{\lambda_{n_i},n_i}$
  converges to $P_n K_{\lambda'}$ in norm.  Hence, $\lim_{i\to\infty}
  \rho(K_{\lambda_{n_i},n_i}) = \rho(K_{\lambda_{n_i},n_i}) = 1$ is an
  eigenvalue of $K_{\lambda'}$ by \thref{krasno}.  Furthermore, a subsequence
  of $(a_{n_i}^{} q_{n_i}^{})$, where $a_n^{} = 1/\|q_n^{}\|_1^{}$, converges
  to an eigenfunction $\tilde{q}$ of $K_{\lambda'}$, and $\tilde{q} \ge 0$
  (but $\tilde{q} \neq 0$).  As there is only one non-negative eigenfunction
  by \thref{existunique}, we conclude $\lambda' = \lambda$ and $\tilde{q} = a
  q$ with $a = 1/\|q\|_1^{}$.  Since this is true for every convergent
  subsequence of $(\lambda_n)$, the claim of part \ref{it:thmreigval} and the
  convergence $a_n^{} q_n^{} \to a q$ follow.
  
  Now, let $n_0$ be sufficiently large such that $\alpha := \inf_{n \ge n_0}
  \lambda_n > 0$.  Then, for $n \ge n_0$,
  \begin{align*}
    \|a_n^{} p_n^{}-(T+\lambda)^{-1} a q\|^{}_1 = 
    \,&\|(P_n^{} T + \lambda_n)^{-1} a_n^{} q_n^{} - 
      (T+\lambda)^{-1} a q\|^{}_1 \\
    \le \,&\|[(P_n^{} T + \lambda_n)^{-1} - (P_n^{} T + \lambda)^{-1}] 
      a_n^{} q_n^{}\|^{}_1 + \\
      &\|(P_n^{} T + \lambda)^{-1} (a_n^{} q_n^{} - a q)\|^{}_1 + \\
      &\|[(P_n^{} T + \lambda)^{-1} - (T+\lambda)^{-1}] a q\|^{}_1 \\
    \le \,&\tfrac{1}{\alpha\lambda} |\lambda-\lambda_n^{}| +
      \tfrac{1}{\lambda} \|a_n^{} q_n^{}-a q\|^{}_1 +
      \tfrac{1}{\lambda^2} \|(\Id - P_n^{}) T a q\|^{}_1 
    \to 0 \es.
  \end{align*}
  With this, $a_n^{} p_n^{} \to a p$ in $\LL{1}(I)$, hence $a_n^{} \to a$ and
  $p_n^{} \to p$, which proves part \ref{it:thmreigvec}.
\end{proof}

\section{Comparison of both methods}
\label{sec:comparison}

Both approaches, the application of the Nystr\"om method in the case
of a compact interval and of the Galerkin method in the case of an
unbounded interval, effectively lead to the same approximation
procedure in our case of the COA model.  First, one chooses
appropriate intervals $I_{n,k}$ and points $t_{n,k} \in I_{n,k}$ (also
for an unbounded interval the use of identical points $t^w_{n,k} =
t^{ux}_{n,k\ell} = t^{uy}_{n,\ell k} = t_{n,k}$ seems reasonable in
many cases).  Then, the operators $T$ and $U$ from \eref{deft} and
\eref{defu}, respectively, are approximated by matrices $\mt{T}_n$ and
$\mt{U}_n$, \cf \eref{deftn}, \eref{defun}, and \eref{rdeftnun}.  For
these, the (finite-dimensional) eigenvalue problem $(\mt{T}_n^{} -
\mt{U}_n^{} + \lambda_n^{}) \vc{p}_n^{} = 0$ is solved.  Here, the
eigenvectors $\vc{p}_n^{}$ are considered as probability
\emph{densities} on $I$.  Then, under the conditions described above,
the eigenvalues $\lambda_n$ converge to $\lambda$ and the measures
corresponding to the $\vc{p}_n^{}$ converge in total variation to the
equilibrium distribution described by the solution $p$ of the original
problem \eref{coaequil}.

The differences between the two approaches lie on the intermediate
technical level of the compact operators $K_\alpha$ and $K_{\alpha,n}$
and the solutions $q$ and $q_n^{}$ of the equivalent eigenvalue
problems \eref{coaequilopk}, \eref{ieigvaleqkn}, and
\eref{reigvaleqkn}.  Here, in the first case we have collectively
compact convergence $K_{\lambda_n,n} \tocc K_\lambda$ going together
with $\|q_n^{} - q\|^{}_\infty \to 0$, whereas in the second case
$\|P_n K_\lambda - K_\lambda\|^{}_Y \to 0$ in $Y = \LL{1}(\RR)$ and
$\|K_{\lambda_n,n} - P_n K_\lambda\|^{}_{Y_n} \to 0$ in the subspaces
$Y_n$ going together with $\|q_n^{} - q\|^{}_1 \to 0$.  On this level,
neither does $\|K_{\lambda_n,n} - K_\lambda\|^{}_\infty \to 0$ hold in
the first case, compare \cite[\thm 12.8]{KreLIE}, nor any kind of
collectively compact convergence in the second.

Both methods may, strictly speaking, only be applied to continuous
mutation kernels $u$.  This excludes, for example,
$\Gamma$-distributions (reflected at the source type), where $u(x,y)
\propto |x-y|^{\Theta-1} \exp(-d\,|x-y|)$, which have poles for $x=y$
if $\Theta \in \ooint{0,1}$ and $d>0$.  These distributions
incorporate biologically desirable properties, such as strong
leptokurticity, and have been used, \eg, in \cite{Hil}.  However,
kernels as the above may be approximated arbitrarily well by
continuous ones in the sense that the norm of the difference
operator---and thus the difference of the largest eigenvalues---gets
arbitrarily small.  Then, the procedures described here may be
applied to these continuous kernels.

\section{Outlook}
\label{sec:outlook}

This article shows that most reasonable COA models can be approximated
arbitrarily well by models with discrete types.  Therefore, one can expect
both model classes to behave quite similarly.  For certain mutation--selection
models with discrete types, a simple maximum principle for the equilibrium
mean fitness $\lambda$ was recently found \cite{HRWB} (see also
\cite{GeBa03,GaGr,BBBK}).  It takes the form
\begin{displaymath}
  \lambda \simeq \sup_{x \in I} \bigl( r(x) - g(x) \bigr)
\end{displaymath}
and holds as an exact identity in a limit of infinitely many types
that densely fill a compact interval $I$.  In the simplest case, a
linear ordering of types is assumed and mutation is taken to only
connect every type $x$ with its two neighbors at rates $u^\pm(x)$.
Then, the function $g$ is given as $g(x) = u^+(x) + u^-(x) - 2
\sqrt{u^+(x) \, u^-(x)}$.  In a subsequent analysis \cite{GaGr}, models
with three types of mutation---and hence six neighbors of every
type---were considered.  For these, $g$ is given as the sum of three
terms of the above pattern (and $x$ has three components), one for
each type of mutation.

In the light of the findings presented here, one may conjecture that
also for certain COA models the above characterization is valid with
an appropriate function $g$.  First steps in \cite{Diss}, both
analytical and numerical, corroborate this conjecture with
\begin{displaymath}
  g(x) = \int_I \bigl( u(x,y) - \sqrt{u(x,y) \, u(y,x)} \bigr) \dd{y} \es,
\end{displaymath}
which generalizes the additive structure of $g$ found in \cite{GaGr}
\wrt a continuum of possible mutations.  The important prerequisite
seems to be the possibility to approximate every local subsystem,
corresponding to a small interval $J \subset I$, by a COA model whose
mutation kernel is of the form $u(x,y) = \exp(\gamma\,(x-y)) \,
h(|x-y|)$.  Then, in a limit $\nu\to\infty$, where $h$ is replaced by
$h_\nu(|x-y|) = \nu \, h(\nu \, |x-y|)$, the above expression seems to
become exact.  A rigorous proof for this statement seems feasible in
the near future.

\begin{ack}
  I thank Michael Baake, Ellen Baake, Tini Garske, Joachim Hermisson, and
  Reinhard B\"urger for discussions and helpful comments on the manu\-script.
  Support by a PhD scholarship of the Studienstiftung des deutschen Volkes is
  gratefully acknowledged.  Furthermore, I express my gratitude to the Erwin
  Schr\"o\-din\-ger International Institute for Mathematical Physics in Vienna
  for support during a stay in December 2002, where the manuscript for this
  article was completed.  Many thanks also go to an anonymous referee for
  hints at how to simplify the main proofs and improve the presentation.
\end{ack}

\let\oldthebibliography\thebibliography
\renewcommand\thebibliography[1]{\oldthebibliography}

\end{document}